\documentclass{article}       
\usepackage{graphicx}
\usepackage{amsmath}

\usepackage{mathtools}
\usepackage{amssymb}

 \usepackage{mathptmx}      
\title{Bivariate rational approximations
	of the general temperature integral
}

\author{
Alireza Aghili, 
              {aghili@iaushiraz.ac.ir, alirezaaghili@yahoo.com}
 \and             
Nadezda Sukhorukova 
 {nsukhorukova@swin.edu.au}
\and
Julien Ugon 
          {julien.ugon@deakin.edu.au}\\
}

\usepackage{hyperref}
\usepackage{titling}
\begin{document}
\maketitle
\begin{abstract}
The non-isothermal analysis of materials with the application of the Arrhenius equation involves temperature integration. If the frequency factor in the Arrhenius equation depends on temperature with a power-law relationship, the integral is known as the general temperature integral. This integral which has no analytical solution is estimated by the approximation functions with different accuracies. In this article, the rational approximations of the integral were obtained based on the minimization of the maximal deviation of bivariate functions. Mathematically, these problems belong to the class of quasiconvex optimization and can be solved using the bisection method. The approximations  obtained in this study are more accurate than all approximates available in the literature.\\ {Keywords: Thermal analysis,  general temperature integral,  bivariate rational 	approximate,  quasiconvex optimization,  bisection method,  uniform 	approximation}\\ {Subclass: 90C26, 90C90, 90C47, 65D15, 65K10, 80A30} \end{abstract}

\section{Introduction}
\label{intro}
Thermal analysis consists of different techniques, which are used for understanding and modeling the thermal behavior of materials. The analysis is often carried out in
non-isothermal conditions with a constant heating rate, and the Arrhenius equation is used for the reaction rate.

The rate of reaction, \(d\alpha/dt\), depends on the temperature and the type of reaction, as given by the following equation:

\begin{equation}\label{eq:1}
	\frac{d\alpha}{dt} = k\left( T \right) f \left( \alpha \right),
\end{equation}
where \(k(T)\) is the rate constant and \(f(\alpha)\) is a function of
conversion. \(k(T)\) can be expressed by the modified Arrhenius equation
\cite{Laidler} as:
\begin{equation}\label{eq:2}
	k(T) = \left( A T^{m} \right)
 \exp\left( - E/RT \right),
\end{equation}
where the term \(\left( A {\ T}^{m} \right)\) is the
temperature-dependent frequency factor, \(A\) and \(m\) are constants,
\(E\) is the activation energy, \(R\) is the universal gas constant and
\(T\) is the absolute temperature, respectively.

The simple form of the Arrhenius equation is a special case with
\(m = 0\), that means the frequency factor is assumed to be constant.
However, for some solid-state reactions \(m\) tends to range from
-1.5 to 2.5 \cite{Criado}.

If \(\beta\) is the constant heating rate, then \(dT = \beta\ dt\), and
it follows that:
\begin{equation}\label{eq:3}
	\frac{d\alpha}{dT} = \frac{A}{\beta} {\ T}^{m}\exp\left( -\frac{E}{RT} \right) f\left( \alpha \right)
\end{equation}

which, upon integration, becomes:

\begin{equation}\label{eq:4}
	F\left( \alpha \right) = \int_{0}^{\alpha}\frac{d\alpha}{f\left( \alpha \right)} = \frac{A}{\beta}\int_{T_{0}}^{T} T^{m}\exp \left(-\frac{E}{RT} \right) {dT} = I_{m} \left( T \right) - I_{m} \left( T_{0} \right),
\end{equation}

where \(T_{0}\) is the initial reaction temperature and
\(I_{m}\left( T \right)\) is:

\begin{equation}\label{ref:5}
	I_{m}\left( T \right)=\frac{A}{\beta}\int_{0}^{T} T^{m} \exp\left( -\frac{E}{RT} \right) dT 
\end{equation}

The temperature integral in equation (\ref{eq:4}) is known as the \emph{general
	temperature integral} \cite{Cai2008}, and if we let \(m = 0\), the well-known
\emph{Arrhenius integral} will be obtained:

\begin{equation}\label{eq:6}
	I_{0} \left( T \right) = \frac{A}{\beta}\int_{0}^{T}\exp \left( -\frac{E}{RT} \right){dT}
\end{equation}

If we let \(x = E/RT\), then:

\begin{equation}\label{eq:7}
	I_{m}\left( T \right) = \frac{A}{\beta}\left( \frac{E}{R} \right)^{m + 1}\int_{x}^{\infty}{\frac{e^{- x}}{x^{m + 2}}dx} = \frac{A}{\beta}\left( \frac{E}{R} \right)^{m + 1}g(m,x)
\end{equation}

and

\begin{equation}\label{eq:8}
	g(m,x) = \int_{x}^{\infty}{\frac{e^{- x}}{x^{m + 2}}dx}
\end{equation}

and for the case of \(m = 0\),

\begin{equation}\label{eq:9}
	g(0,x) = \int_{x}^{\infty}{\frac{e^{- x}}{x^{2}}dx}
\end{equation}

The integrals, in general form, have no analytical solution. Integration
by parts results in \cite{Criado}:

\begin{equation}\label{eq:10}
	g(m,x) = 
	\frac{e^{- x}}{x^{m + 2}} 
	\left(1 - \frac{\left( m + 2 \right)}{x} + \frac{\left( m + 2 \right)\left( m + 3 \right)}{x^{2}}
	- 	\frac{\left( m + 2 \right)\left( m + 3 \right)\left( m + 4 \right)}{x^{3}}
	 + -\ldots 
	\right) 
\end{equation}

The series expansion shows infinite terms, which means \(g(m,x)\) is not
an elementary function. Usually different elementary approximation
functions, with finite terms and different accuracies, are employed to
evaluate the integrals. In the literature, \(g(0,x)\) has been
substituted by several series, asymptotic, rational, or semi-empirical
functions. Deng et al \cite{Deng}, Flynn and Wall \cite{Flynn}, Lyon \cite{Lyon}, Pérez-Maqueda et al \cite{Perez} and Órfão \cite{Orfao}  have provided excellent
reviews of these functions. The function \(g(m,x)\) has also been
estimated with different series, rational or semi-empirical functions \cite{Cai2007-1,Cai2008,Cai2007-2,Capela,Casal,Chen2007,Chen2009-1,Chen2009-2,Gorbachev,Lei,Wanjun2009,Wanjun2005,Xia}  
with different accuracies.

In this study, new rational approximations of the general temperature
integral are obtained based on the minimization of the maximal deviation
of bivariate functions. Minimization of maximal deviation is also known
as uniform (or Chebyshev) approximation.

The theory of univariate uniform approximation by rational function was
originally developed by Achieser \cite{Achieser}. Then the results were
extended by Loeb \cite{Loeb} to the case when the approximations are ratios of linear forms and  the elementary functions of these linear forms  are not
limited to monomials. This is also known as generalized rational
approximation. The procedure was extended by Diaz Millán et al \cite{Millan}
 to multivariate approximation.

Rational approximations are viewed as a natural extension to polynomial
approximations \cite{Chebyshev,Nurnberger,Remez} (just to name a few). Rational approximations are more flexible than polynomials and therefore are suitable for
approximating non-smooth and non-Lipschitz functions (functions with
abrupt changes).

The numerical methods for rational and generalized rational
approximation can be divided into two groups. The first group aims at
finding "nearly optimal" solutions and their constructions are based
on Chebyshev polynomials (Trefethen \cite{Trefethen}). These methods are fast
and efficient, but the solutions are not always optimal, especially in
the case of multivariate approximation. The second group of methods are
searching for optimal solutions \cite{Loeb,Millan}. This research is based on
the computational procedure, developed by Diaz Millán et al \cite{Millan}.
This research targets multivariate rational approximations, while the
procedure from \cite{Millan} is also suitable for generalized rational
approximations. The results of numerical experiments demonstrate that,
despite its simplicity, this procedure is robust and accurate. The
procedure is described in the next section.

\section{Theory}\label{theory}

The rational approximates, in this study, have been obtained to estimate the function \(h(m,x)\),
the term in the last brace in equation~(\ref{eq:10}):

\begin{equation}\label{eq:12}
	h(m,x) =  {e^{ x}}{x^{m + 2}} g(m,x)
\end{equation}

It was mentioned that for some solid-state reactions \(m\) varies
between \(- 1.5\) and \(2.5\) \cite{Criado}, however, for pyrolytic reactions or gas combustion reactions, a wider range, between \(- 4\) and
\(3.2\), has also been reported \cite{Kovacs,Turanyi}. Therefore, we can assume that \(m\) varies generally between \(- 4\) and \(4\). On the other hand, Starink \cite{Starink2003} showed that the majority of reactions occur in the range \(15 < x < 60\) and only the range \(9 < x < 100\) is of practical significance for approximation functions. But some reactions occur at \(x\) as low as \(4.6\) \cite{Starink2018}. Therefore, the intervals \(- 4 \leq m \leq 4\) and \(4 \leq x \leq 100\) were used to find the approximation functions.

It was demonstrated in \cite{Millan} that the objective functions (i.e., the
functions that are subject to minimization) in the corresponding
optimization problems are quasiconvex. The notion of quasiconvexity was
originally introduced by Bruno de Finetti's \cite{deFinetti}, a mathematician
who was working with the area of financial mathematics. The term,
however, was introduced later as a possible extension to convexity.

For
a given value \(u\) the set \(S_{u} = \left\{ X:f(X) \leq u \right\}\)
is called a sublevel set of the function $f$.
A quasiconvex function is a function whose sublevel sets are all convex.

The method of bisection for quasiconvex functions \cite{Boyd} is an algorithm for minimizing a quasiconvex function subject to convex (or linear) constraints or without constraints. The procedure terminates in a finite
number of steps and the solution is within a given precision of the optimal solution. The main difficulty of this procedure is the solution of the so-called convex feasibility problem, which may be a difficult
problem \cite{Bauschke}. Diaz Millán et al \cite{Millan} demonstrated that in the
case of rational and generalized rational approximation on finite grids, the constraints are linear and the convex feasibility problems can be reduced to solving linear programming problems and therefore can be
efficiently implemented and the solution is within a prespecified
precision of the optimal solution.

The problem of approximating the function \(h\left (m,x\right )\) with a rational function \(r\left (m,x\right )\), where \(r = p/q\), the numerator is a polynomial \(p\left (m,x\right )\) of degree
\(n_{1}\): \(p \in \Pi_{n_{1}}\) and the denominator \(q\left (m,x\right )\) is a polynomial of degree \(n_{2}\): \(q \in \Pi_{n_{2}}\) can be formulated as follows:

\begin{equation*}\label{eq:obj}
{\rm minimize}~ u
\end{equation*}
\ subject\ to\ 
\begin{equation*}\label{eq:constraints}
\left| h\left( m,x \right) - \frac{p\left( m,x \right)}{q\left( m,x \right)} \right| \leq u,\ \forall\ m\  \in \ \lbrack -4,4\rbrack,x \in \ \lbrack 4,100\rbrack\ 
\end{equation*}

The value of \(u\) represents the maximal deviation. For a fixed \(u\),
the set of feasible points (that is, the \emph{sublevel set}) satisfying
the constraint:

\begin{equation}\label{eq:11}
	S_{u} = \left\{ m,x:\left| h\left( m,x \right) - \frac{p\left( m,x \right)}{q\left( m,x \right)} \right| \leq u \right\}
\end{equation}

is convex, and so the problem is quasiconvex. Therefore, the aim is to
find the smallest \(u\) such that \(S_{u}\) is nonempty. This leads to a
natural application of the bisection algorithm as follows:

\begin{enumerate}
	\item
	Set \(u_{-} = 0\) and \(u_{+} = \parallel h \parallel_{\infty}\).
	\item
	For \(u = \frac{u_{-} + u_{+}}{2}\), solve the convex feasibility
	problem \cite{Millan}. If this
	problem has a solution, the sublevel set \(S_{u}\) is nonempty, and so
	set \(u_{+} = u\). If the problem has no solution, the set \(S_{u}\)
	is empty, and so set \(u_{-} = u\).
	\item
	Repeat Step 2 until \(u_{+} - u_{-}\) is below a desired accuracy.
\end{enumerate}

At every step of the algorithm, it is easy to see that the solution to
the problem \(u^{*} \in \left\lbrack u_{-},u_{+} \right\rbrack\), and
that the size of that interval is halved at every iteration. The
algorithm returns an approximate solution to the problem of
approximating the function \(h\) with a rational function \(r\).

In this paper we are using a discretization of the intervals, and
therefore the set \(S_{u}\) is a polytope, delimited by a finite number
of hyperplanes (linear equations). The convex feasibility subproblem
solved at each iteration of the bisection algorithm can be solved with
any linear programming toolbox.

\section{Results and discussion}\label{results-and-discussion}

The rational approximations of
\(h(m,x)\) with the same degrees of polynomials in numerators and denominators \(\left (n_{1}=n_{2} =n \right)  \)  were calculated by using the
procedure explained in the previous section. Therefore, the approximates to
\(g(m,x)\), denoted here by \(g_{n}(m,x)\) will have the following
general form:

\begin{equation}\label{eq:13}
	g_{n}(m,x) = \left( \frac{e^{- x}}{x^{m + 2}} \right)\left( \frac{\sum_{i = 0}^{n}{\sum_{j = 0}^{n - i}{a_{ij}x^{i}m^{j}}}}{\sum_{i = 0}^{n}{\sum_{j = 0}^{n - i}{b_{ij}x^{i}m^{j}}}} \right)
\end{equation}

The approximation functions of the general temperature integral with different degrees \(n = 1\sim 4\) are listed in table \ref{table01}
and the corresponding coefficients of the polynomials in the numerators (\(a_{ij}\))
and denominators (\(b_{ij}\)), are shown in tables \ref{table 1} to \ref{table 4}.

\begin{table}
	\caption{New approximates to the general temperature integral with different degrees $(n)$}\label{table01}
	\renewcommand{\arraystretch}{1.9}
	\begin{tabular}{|c|l|}
		\hline
		$n$ & $g_{n}(m,x)$\\
		\hline
		$1$ & $\left( \frac{e^{- x}}{x^{m + 2}} \right)\left(\frac{a_{00}+a_{10}x+a_{01}m}{b_{00}+b_{10}x+b_{01}m} \right)$\\
		\hline
		$2$ & $\left( \frac{e^{- x}}{x^{m + 2}} \right)\left(\frac{a_{00}+a_{10}x+a_{01}m+a_{20}x^{2}+a_{02}m^{2}+a_{11}xm}{b_{00}+b_{10}x+b_{01}m+b_{20}x^{2}+b_{02}m^{2}+b_{11}xm} \right)$\\
		\hline
		$3$ & $\left( \frac{e^{- x}}{x^{m + 2}} \right)\left(\frac{a_{00}+a_{10}x+a_{01}m+a_{20}x^{2}+a_{02}m^{2}+a_{11}xm+a_{30}x^{3}+a_{03}m^{3}+a_{21}x^{2}m+a_{12}xm^{2}}{b_{00}+b_{10}x+b_{01}m+b_{20}x^{2}+b_{02}m^{2}+b_{11}xm+b_{30}x^{3}+b_{03}m^{3}+b_{21}x^{2}m+b_{12}xm^{2}} \right)$\\
		\hline
		$4$ & $\left( \frac{e^{- x}}{x^{m + 2}} \right)\left(\frac{a_{00}+a_{10}x+a_{01}m+a_{20}x^{2}+a_{02}m^{2}+a_{11}xm+a_{30}x^{3}+a_{03}m^{3}+a_{21}x^{2}m+a_{12}xm^{2}+a_{40}x^{4}+a_{04}m^{4}+a_{31}x^{3}m+a_{13}xm^{3}+a_{22}x^{2}m^{2}}{b_{00}+b_{10}x+b_{01}m+b_{20}x^{2}+b_{02}m^{2}+b_{11}xm+b_{30}x^{3}+b_{03}m^{3}+b_{21}x^{2}m+b_{12}xm^{2}+b_{40}x^{4}+b_{04}m^{4}+b_{31}x^{3}m+b_{13}xm^{3}+b_{22}x^{2}m^{2}} \right)$\\
		\hline
	\end{tabular}
\end{table}

\begin{table}
	\caption{The coefficients of polynomials in \(g_{1}(m,x)\)}
	\label{table 1}
	\begin{tabular}{|c|c|c|c|}
		\hline
		$a_{00}$ & 0.237276056849810 & $b_{00}$ & 1 \\
		\hline
		$a_{10}$ & 0.388591025647952 & $b_{10}$ & 0.386946448530584 \\
		\hline
		$a_{01}$ & -0.039895879080345 & $b_{01}$ & 0.337632268754683 \\
		\hline
	\end{tabular}	
\end{table}

\begin{table}
	\caption{The coefficients of polynomials in \(g_{2}(m,x)\)}
	\label{table 2}
	\begin{tabular}{|c|c|c|c|}
		\hline
		$a_{00}$ & 0.091419664846862 & $b_{00}$ & 1 \\
		\hline
		$a_{10}$ & 0.592563087097409 & $b_{10}$ & 0.875127841491978 \\
		\hline
		$a_{01}$ & -0.024580517356591 & $b_{01}$ & 0.610076977863149 \\
		\hline
		$a_{20}$ & 0.141353166023985 & $b_{20}$ & 0.141350626515119 \\
		\hline
		$a_{02}$ & 0.001839998211580 & $b_{02}$ & 0.092119571681669 \\
		\hline
		$a_{11}$ & 0.089269674860906 & $b_{11}$ & 0.230585361263812 \\
		\hline
	\end{tabular}	
\end{table}

\begin{table}
	\caption{The coefficients of polynomials in \(g_{3}(m,x)\)}
	\label{table 3}
	\begin{tabular}{|c|c|c|c|}
		\hline
		$a_{00}$ & 0.014728216211734 & $b_{00}$ & 1 \\
		\hline
		$a_{10}$ & 0.782330940685156 & $b_{10}$ & 1.838698460838684 \\
		\hline
		$a_{01}$ & -0.014062001771903 & $b_{01}$ & 0.929983297897110 \\
		\hline
		$a_{20}$ & 0.609319528759108 & $b_{20}$ & 0.771894891782121 \\
		\hline
		$a_{02}$ & 0.002525329133302 & $b_{02}$ & 0.287579275733519 \\
		\hline
		$a_{11}$ & 0.246540819422722 & $b_{11}$ & 0.974064696570030 \\
		\hline
		$a_{30}$ & 0.081289506456857 & $b_{30}$ & 0.081289480719039 \\
		\hline
		$a_{03}$ & -0.000136443018005 & $b_{03}$ & 0.029552723377583 \\
		\hline
		$a_{21}$ & 0.099632579649176 & $b_{21}$ & 0.180921556157630 \\
		\hline
		$a_{12}$ & 0.029904486224604 & $b_{12}$ & 0.129584150595668 \\
		\hline
	\end{tabular}	
\end{table}

\begin{table}
	\caption{The coefficients of polynomials in \(g_{4}(m,x)\)} 
	\label{table 4} 
	\begin{tabular}{|c|c|c|c|}
		\hline
		$a_{00}$ & -0.122782923809827 & $b_{00}$ & 1E-05 \\
		\hline
		$a_{10}$ & 0.354996648791790 & $b_{10}$ & 2.235445022592591 \\
		\hline
		$a_{01}$ & -0.016734436667053 & $b_{01}$ & -0.209413742177965 \\
		\hline
		$a_{20}$ & 1.189400695088519 & $b_{20}$ & 1.765382336034014 \\
		\hline
		$a_{02}$ & 0.011686328539599 & $b_{02}$ & -0.279918727330068 \\
		\hline
		$a_{11}$ & -0.268306511298717 & $b_{11}$ & 0.656323656760673 \\
		\hline
		$a_{30}$ & 0.300563881809569 & $b_{30}$ & 0.325660976846066 \\
		\hline
		$a_{03}$ & -0.001568905108973 & $b_{03}$ & -0.111746651247844 \\
		\hline
		$a_{21}$ & -0.034919185591078 & $b_{21}$ & 0.224710088573674 \\
		\hline
		$a_{12}$ & -0.090995975658326 & $b_{12}$ & -0.170181138653749 \\
		\hline
		$a_{40}$ & 0.012548386262555 & $b_{40}$ & 0.012548392422441 \\
		\hline
		$a_{04}$ & 0.000071049885199 & $b_{04}$ & -0.013889020864874 \\
		\hline
		$a_{31}$ & -0.014175730850403 & $b_{31}$ & -0.001627152422281 \\
		\hline
		$a_{13}$ & -0.014071092519943 & $b_{13}$ & -0.049775575854942 \\
		\hline
		$a_{22}$ & -0.035706031933476 & $b_{22}$ & -0.049888788836345 \\
		\hline
	\end{tabular}	
\end{table}

\subsection{Accuracy of the rational approximation functions}

To evaluate the accuracy of the obtained rational approximations, they
should be compared with the high precision values of the integral. The
relative deviation is defined by the following equation:

\begin{equation}\label{eq:14}
	\epsilon = \left( \frac{g_{n}(m,x)\ }{g(m,x)} - 1 \right)
\end{equation}

In the above equation, \(g(m,x)\) must be computed with high precision
by numerical integration. The intervals \(- 4 \leq m \leq 4\) and
\(4 \leq x \leq 100\) with the step sizes of \(\Delta x = 1\) and \(\Delta m = 0.1\) were used for evaluation of the accuracy of the
approximation functions.

The accuracy of the rational approximations to the general temperature
integral (equations in table ~\ref{table01}), was evaluated and the results are
shown in figures S1 to S4 in the supplementary information file. In
these figures, the relative error has been illustrated as a function of
\(m\) and \(x\). The maximum values of the relative deviations,
\(\left| \epsilon \right|_{\max}\), and sum of squared errors, SSE, are
also listed in table~\ref{table 5}. The rational approximations show very high
accuracies, particularly when higher degrees are employed.

\begin{table}
	\caption{Accuracy of new approximation functions}
	\label{table 5}
	\begin{tabular}{|c|c|c|}
		\hline
		Degree & $\vert \epsilon \vert _{\max }$	& SSE\\
		\hline
		$n=1$ & 1.12E-02 & 1.34E-01 \\
		\hline
		$n=2$ & 6.26E-05 & 3.21E-06 \\
		\hline
		$n=3$ & 1.72E-06 & 1.80E-09 \\
		\hline
		$n=4$ & 6.18E-07 & 3.77E-10 \\
		\hline
	\end{tabular}
\end{table}

\subsection{Application and comparison with other approximation
	functions}\label{application-and-}

The Arrhenius integral has been widely used for kinetic modeling of
thermal analysis of materials. For example, in a single step reaction
with relatively constant activation energy (\(E\)), to find the best
model for the thermal behavior of a substance, the functions
\(y(\alpha)\) and \(z(\alpha)\) are evaluated by the following equations
and compared with master plots provided in the literature \cite{Vyazovkin2011}.

\begin{equation}\label{eq:15}
	y\left( \alpha \right) = Af\left( \alpha \right) = \beta\left( \frac{d\alpha}{dT} \right)_{\alpha}\exp\left( \frac{E}{RT_{\alpha}} \right)
\end{equation}

\begin{equation}\label{eq:16}
	z\left( \alpha \right) = f\left( \alpha \right)F\left( \alpha \right) = \frac{A}{\beta}\frac{E}{R}f\left( \alpha \right)g(0,x)
\end{equation}

The Arrhenius integral, \(g(0,x)\), in this equation is usually
estimated by the approximation functions.

Another example for the application of the Arrhenius integral is the
Vyazovkin's advanced isoconversional method \cite{Vyazovkin1996} in which the
activation energy is determined as a function of conversion. For a
series of runs performed at different heating rates, the \(E_{\alpha}\)
value can be determined by minimizing the following function:

\begin{equation}\label{eq:17}
	\Phi(E_{\alpha}) = \sum_{i = 1}^{n}{\sum_{j \neq i}^{n}\left\lbrack \frac{\beta_{j}I\left( E_{\alpha},T_{\alpha,i} \right)}{\beta_{i}I\left( E_{\alpha},T_{\alpha,j} \right)} \right\rbrack},
\end{equation}

where

\begin{equation}\label{eq:18}
	I\left( E_{\alpha},T_{\alpha} \right) = \int_{T_{\alpha-\Delta\alpha}}^{T_{\alpha}}\exp\left( -\frac{E_{\alpha}}{RT} \right){dT} = \frac{E_{\alpha}}{R}[g(0,x_{\alpha})-g(0,x_{\alpha-\Delta\alpha})]
\end{equation}

In this method, the univariate function \(g(0,x)\) is evaluated numerically or estimated by
the approximation functions \cite{Galukhin,Vyazovkin1996}.

In this section, a comparison has been made between the highly accurate
approximations of \(g(0,x)\) proposed in the literature and new approximates
obtained in this study. Some highly accurate approximates to \(g(0,x)\)
are listed in table 6. As mentioned earlier, excellent reviews of the
approximates to the Arrhenius integral have been provided in \cite{Deng,Flynn,Lyon,Perez,Orfao}. Órfão \cite{Orfao} used the O model as the highest accurate model among the approximates.  Deng et al \cite{Deng} showed that The O and J models are the most accurate approximates to the Arrhenius integral. Therefore, it would be sufficient to compare new approximates with only the O and J models. On the other hand, the SY model has been used more extensively by the
researchers.

\begin{table}
	\caption{Approximates to the Arrhenius integral}
	\label{table 6}
	\renewcommand{\arraystretch}{1.8}
	\begin{tabular}{|l|c|l|}
		\hline
		Reference & Model & Approximation function\\
		\hline
		Ji \cite{Ji} & J & $
		(\frac{e^{-x}}{x^{2}})(\frac{x^{2}+16.99864x+3.65517\ln 
			x+5.41337}{x^{2}+18.99977x+3.43593\ln x+38.49858})$ \\
		\hline
		Órfão \cite{Orfao} (4$^{th}$ degree) & O & $
		(\frac{e^{-x}}{x^{2}})(\frac{0.9999936x^{4}+7.5739391x^{3}+12.4648922x^{2}+3.6907232x}{x^{4}+9.5733223x^{3}+25.6329561x^{2}+21.0996531x+3.9584969})
		$ \\
		\hline
		Senum \& Yang \cite{Senum} (4$^{th}$ degree) & SY$^{*}$ & $
		(\frac{e^{-x}}{x^{2}})(\frac{x^{4}+18x^{3}+86x^{2}+96x}{x^{4}+20x^{3}+120x^{2}+240x+120})
		$ \\
		\hline
	\end{tabular}
	\textsuperscript{*} Some authors, including Senum and Yang \cite{Senum} have erroneously quoted the coefficient of $x^2$ in numerator as 88 and not the correct value of 86 \cite{Starink2003}.
\end{table}

The approximation functions (table~\ref{table 6}) and the new approximates proposed
in this study, (tables \ref{table01} to \ref{table 4} with \(m = 0\)) were compared in the
interval \(4 \leq x \leq 100\) with a step size of \(\Delta x = 1\), and the
results of the comparison of these models, including the SSE and the
maximum relative deviation \(\left| \epsilon \right|_{\max}\), are shown
in table~\ref{table 7}. The results reveal that the SSE and
\(\left| \epsilon \right|_{\max}\) of the new approximates with degrees
3 and 4 are less than those of the J, O and SY models. Figure S5 in the
supplementary information file shows comparison of the relative
deviations of these models. The J and SY models have less deviations
than new approximates only for \(x > 18\) and \(x > 20\), respectively,
but if the whole interval \(4 \leq x \leq 100\) is being considered, the
new approximates, \(g_{3}\left( m,x \right)\) and
\(g_{4}\left( m,x \right)\), are more accurate.

\begin{table}
	\caption{Comparison of different models for evaluation of the Arrhenius
		integral}
	\label{table 7}
	\begin{tabular}{|c|c|c|c|}
		\hline
		Model Name & Equation in: & SSE & $\vert \epsilon \vert _{\max }$\\
		\hline
		J & Table \ref{table 6} & 3.45E-11 & 5.66E-06 \\
		\hline
		O & Table \ref{table 6} & 7.25E-11 & 1.87E-06 \\
		\hline
		SY & Table \ref{table 6} & 7.86E-09 & 8.15E-05 \\
		\hline
		$g_{1}(m,x)$ & Table \ref{table01} & 1.89E-03 & 6.79E-03 \\
		\hline
		$g_{2}(m,x)$ & Table \ref{table01} & 4.62E-08 & 3.95E-05 \\
		\hline
		$g_{3}(m,x)$ & Table \ref{table01} & 1.33E-11 & 8.89E-07 \\
		\hline
		$g_{4}(m,x)$ & Table \ref{table01} & 6.72E-12 & 3.95E-07 \\
		\hline
	\end{tabular}
\end{table}

The general temperature integral has also been used in the literature
for thermal analysis of the materials. For example, an integral method
of combined kinetic analysis of reactions (ICKA) has been proposed by
Casal and Marbán \cite{Casal}. In this method, \(f(\alpha)\) is estimated by a
general relationship as follows:

\begin{equation}\label{eq:19}
	f\left( \alpha \right) = Z\left\lbrack 1 - \left( 1 - \alpha^{a} \right)^{b} \right\rbrack^{1 - c}\left( 1 - \alpha^{a} \right)^{1 - b}\alpha^{1 - a}
\end{equation}

in which \(Z\), \(a\), \(b\), and \(c\) are parameters that depend on
the kinetic model. Then the conversion (\(\alpha\)) will be calculated
by the following relationship:

\begin{equation}\label{eq:20}
	\alpha = \left\{ 1 - \left\lbrack 1 - \left( \frac{abcZA}{\beta}\left( \frac{E}{R} \right)^{m + 1}g(m,x) \right)^{1/c} \right\rbrack^{1/b} \right\}^{1/a}
\end{equation}

In this equation, \(g(m,x)\) is the general temperature integral that
is estimated by the approximation functions.

A comparison has also been made between the new approximations with the
approximations proposed in the literature. Table~\ref{table 8} shows the
approximation functions in the literature. One should take into account
that for the Ch1 model \cite{Chen2007} and X model \cite{Xia}, the fourth-degree
rational functions are used here, and in the X model \cite{Xia} the
coefficients of the polynomial in numerator are explicitly reported only for limited values of \(m\), indicated in table~\ref{table 9}, and the evaluation is based on these values. As we know, the approximations were obtained for the range of
\(- 1.5 \leq m \leq 2.5\). Therefore, two cases were considered for
comparing the approximation functions; the first in the interval
\(- 1.5 \leq m \leq 2.5\) and the second in the interval
\(- 4 \leq m \leq 4\), and for both cases, the range of
\(4 \leq x \leq 100\) with the step sizes of \(\Delta x = 1\) and \(\Delta m = 0.1\) were used. The figures S6 to S19 in the supplementary information file,
show the deviations of the approximation functions for the first case,
\(- 1.5 \leq m \leq 2.5\) and \(4 \leq x \leq 100\).

For both cases, the maximum deviations,
\(\left| \epsilon \right|_{\max}\), and SSE were obtained for each model
and the results are listed in table~\ref{table 10}. The results show that except for
the first degree, \(g_{1}(m,x)\) , the new approximates obtained in this
study (\(g_{n}(m,x)\) with \(n = 2\),\(3\), and \(4\)) have less values
of SSE and \(\left| \epsilon \right|_{\max}\) compared to the other
approximates. Comparison of the figures S1 to S4, and S6 to S19, in the
supplementary information file, confirms that the new approximates have
less maximal deviations. Therefore, new approximation functions are more
accurate than the approximates proposed in the literature.

\begin{table}
	\caption{Approximation functions of the general temperature integral}
	\label{table 8}
	\renewcommand{\arraystretch}{1.8}
	\begin{tabular}{|l|c|l|}
		\hline
		Reference & Model name & Approximation function \\
		\hline
		Gorbachev \cite{Gorbachev} & G & $
		(\frac{e^{-x}}{x^{m+2}})(\frac{1}{(1+\frac{m+2}{x})})$ \\
		\hline
		Wanjun \cite{Wanjun2005} & W1 & $(\frac{e^{-x}}{x^{m+2}})(\frac{1}{\lbrace 
			1+(m+2)(0.00099441+\frac{0.93695599}{x})\rbrace })$ \\
		\hline
		Wanjun \cite{Wanjun2009} & W2 & $
		\frac{exp(-0.18887(m+2)-(1.00145+0.00069m)x)}{x^{0.94733(m+2)}}$ \\
		\hline
		Cai \cite{Cai2007-1} & C1 & $(\frac{e^{-x}}{x^{m+2}})(\frac{0.99954 x-0.044967 
			m+0.58058}{x+0.94057 m+2.5400 })$ \\
		\hline
		Cai \cite{Cai2007-2} & C2 & $(\frac{e^{-x}}{x^{m+2}})(\frac{1.0002486 
			x+0.2228027\ln x-0.05241956 m+0.2975711}{x+0.2333376\ln x+0.9496628 
			m+2.2781591 })$ \\
		\hline
		Cai \cite{Cai2008} & C3 & $(\frac{e^{-x}}{x^{m+2}})(\frac{x-0.054182 
			m+0.65061}{x+0.93544 m+2.62993 })$ \\
		\hline
		Chen \cite{Chen2007} & Ch1 & $
		(\frac{e^{-x}}{x^{m+2}})(\frac{x^{4}+3(m+2)x^{3}+(3m+1)(m+2)x^{2}+m(m-1)(m+2)x}{x^{4}+4(m+2)x^{3}+6(m+1)(m+2)x^{2}+4m(m+1)(m+2)x+(m-1)(m)(m+1)(m+2)})
		$ \\
		\hline
		Chen \cite{Chen2009-1} & Ch2 & $
		\frac{exp(-(0.16656m+0.39329)-(1.00147+0.00057m)x)}{x^{1.89021+0.95479m}}
		$ \\
		\hline
		Chen \cite{Chen2009-2} & Ch3 & $
		(\frac{e^{-x}}{x^{m+2}})(\frac{x}{(1.00141+0.0006m)x+(1.89376+0.95276m) 
		})$ \\
		\hline
		Chen \cite{Chen2009-2} & Ch4 & $
		(\frac{e^{-x}}{x^{m+2}})(\frac{x+(0.74981-0.06396m)}{(1.00017+0.00013m)x+(2.73166+0.92246m) 
		})$ \\
		\hline
		Capela \cite{Capela} & Cp & $
		(\frac{e^{-x}}{x^{m+2}})(\frac{(2-\sqrt{2})}{4}(\frac{x}{x+2+\sqrt{2}})^{m+2}+\frac{(2+\sqrt{2})}{4}(\frac{x}{x+2-\sqrt{2}})^{m+2})
		$ \\
		\hline
		Xia \cite{Xia} & X & $
		(\frac{e^{-x}}{x^{m+2}})(\frac{x^{4}+a_{3}x^{3}+a_{2}x^{2}+a_{1}x}{x^{4}+16x^{3}+72x^{2}+96x+24})
		$ \\
		\hline
		Casal \cite{Casal} & Cs & $(\frac{e^{-x}}{x^{m+2}})(\frac{ x-0.05924479 
			m+0.62385968}{x+0.92755595 m+2.59746116 })$ \\
		\hline
		Lei \cite{Lei} & L & $
		(\frac{e^{-x}}{x^{m+2}})(\frac{\sqrt{(x+m+1)^{2}+4x}-(x+m+1)}{2})$ \\
		\hline
	\end{tabular}
\end{table}

\begin{table}
	\caption{The coefficients of X model \cite{Xia}}
	\label{table 9}
	\begin{tabular}{|c|c|c|c|}
		\hline
		$m$ & $a_{3}$ & $a_{2}$ & $a_{1}$ \\
		\hline
		-1 & 15 & 58 & 50 \\
		\hline
		-0.5 & 14.5 & 51.75 & 34.875 \\
		\hline
		0 & 14 & 46 & 24 \\
		\hline
		0.5 & 13.5 & 40.75 & 16.625 \\
		\hline
		1 & 13 & 36 & 12 \\
		\hline
		2 & 12 & 28 & 8 \\
		\hline
	\end{tabular}
\end{table}

\begin{table}
	\caption{Comparison of different models for evaluation of the general
		temperature integral}
	\label{table 10}
	\begin{tabular}{|c|c|c|c|c|}
		\hline
		& \multicolumn{2}{c|}{$-1.5\leq m \leq 2.5$} &  \multicolumn{2}{c|}{$-4 \leq m \leq 4$}   \\
		\hline
		Model Name& SSE & $\vert \epsilon \vert _{\max }$ & SSE & $\vert \epsilon \vert _{\max }$ \\
		\hline
		G & 2.67E-01 & 5.58E-02 & 8.45E-01 & 2.31E-01 \\
		\hline
		W1 & 5.11E-02 & 2.79E-02 & 2.49E-01 & 1.62E-01 \\
		\hline
		W2 & 1.81E+00 & 1.76E-01 & 3.61E+00 & 2.20E-01 \\
		\hline
		C1 & 1.11E-03 & 7.89E-03 & 2.34E-02 & 5.42E-02 \\
		\hline
		C2 & 3.53E-05 & 1.21E-03 & 2.07E-02 & 5.76E-02 \\
		\hline
		C3 & 2.29E-03 & 1.02E-02 & 1.24E-02 & 3.71E-02 \\
		\hline
		Ch1 & 2.31E-04 & 3.36E-03 & 3.06E-01 & 2.97E-01 \\
		\hline
		Ch2 & 1.82E+00 & 2.02E-01 & 4.90E+00 & 2.65E-01 \\
		\hline
		Ch3 & 7.16E-02 & 3.26E-02 & 3.51E-01 & 1.84E-01 \\
		\hline
		Ch4 & 1.03E-02 & 1.88E-02 & 1.93E-02 & 1.91E-02 \\
		\hline
		Cp & 4.39E-02 & 5.58E-02 & 2.12E-01 & 1.04E-01 \\
		\hline
		X $^{*}$ & $-$ & 6.14E-04 & $-$ & $-$ \\
		\hline
		Cs & 1.50E-03 & 7.21E-03 & 1.02E-02 & 3.60E-02 \\
		\hline
		L & 1.36E-03 & 4.43E-03 & 1.21E-02 & 3.90E-02 \\
		\hline
		$g_{1}(m,x)$ & 7.67E-02 & 7.13E-03 & 1.34E-01 & 1.12E-02 \\
		\hline
		$g_{2}(m,x)$ & 1.73E-06 & 6.25E-05 & 3.21E-06 & 6.26E-05 \\
		\hline
		$g_{3}(m,x)$ & 6.80E-10 & 1.29E-06 & 1.80E-09 & 1.72E-06 \\
		\hline
		$g_{4}(m,x)$ & 2.13E-10 & 5.19E-07 & 3.77E-10 & 6.18E-07 \\
		\hline
	\end{tabular}\\
	\textsuperscript{*} The value of \(\left| \epsilon \right|_{\max}\) for
	X model corresponds to the values of \(m\) indicated in table~\ref{table 9}.
\end{table}

As a consequence, it has been shown that the rational approximates to the general temperature integral obtained by minimization of the maximal deviation of bivariate functions are more accurate than all approximates available in the literature. We hope that the new approximates would help the researchers
to develop non-isothermal kinetic analysis of materials.

\section{Conclusion}\label{conclusion}

Application of the Arrhenius equation in the non-isothermal analysis
with linear heating involves the calculation of a temperature integral.
This integral is known as the Arrhenius integral if the frequency factor
in the Arrhenius equation is constant, or the general temperature
integral if the frequency factor depends on temperature. The 
 integrals have no analytical solutions and the researchers
often use the approximation functions to estimate them. In this article,
by using the minimization of the maximal deviation of bivariate
functions the rational approximations of the general temperature
integral with excellent accuracies were obtained and it was shown that new approximations are more accurate that all approximates in the literature. 

 

\subsection*{Funding} 
This research was supported by the Australian Research Council (ARC),  Solving hard Chebyshev approximation problems through nonsmooth analysis (Discovery Project DP180100602).








\end{document}